\documentclass[reqno]{amsart}
\usepackage{amsmath, amssymb, amscd, eepic, conformal, mystyle}

\begin{document}
\bibliographystyle{plain}

\author{Michael Roitman}


\title{On free conformal and vertex algebras}

\date\today


\thanks{Partially supported by NSF grant DMS-9704132}


\maketitle


Vertex algebras and conformal algebras have recently attracted 
a lot of attention due to their connections with physics and Moonshine
representations of the Monster. See, for example,  \cite{bor}, \cite{flm}, 
\cite{kac_loc}, \cite{kac2}, \cite{lz_moon}.

In this paper we describe bases of free conformal and free vertex
algebras (as introduced in \cite{bor}, see also \cite{lz_preprint}).

All linear spaces are over a field $\Bbbk$ of characteristic 0.
Throughout this paper $\Z_+$ will stand for the set of non-negative
integers.

In \S1 and \S2 we give a review of conformal and vertex algebra
theory. All statements is these sections are either in \cite{fhl},
\cite{kac_loc}, \cite{kac_fd}, \cite{kac2}, \cite{lz_cqoa}, 
\cite{lz_preprint} or easily follow from results therein. In \S3 we
investigate free conformal and vertex algebras. 

\bigskip
\section{Conformal algebras}
\subsection{Definition of conformal  algebras}
\label{sec:def} 
We first recall some basic definitions and constructions, see
 \cite{kac_fd}, \cite{kac2}, \cite{lz_cqoa}, \cite{lz_preprint}. 
The main object of investigation
is defined as follows:

\begin{Def}\label{dfn:thedef}
A {\it Conformal algebra} is a linear space $C$ endowed with a linear operator
$D:C \to C$ and a sequence of bilinear products 
$\encirc n: C\otimes C \to C, \ n \in \Z_+$, such that
for any $a,b \in C$ one has
\begin{itemize}
\item[(i)](locality) There is a non-negative integer $N = N(a,b)$
such that $a \encirc n b = 0$ for any $n \ge N$;
\item[(ii)]
$D(a \encirc n b) = (Da)\encirc n b + a \encirc n (Db)$;
\item[(iii)]
$(Da)\encirc n b = -n a\encirc{n-1} b$.
\end{itemize}
\end{Def}

\bigskip
\subsection{Spaces of power series} 
\label{sec:series}
Now let us discuss the main motivation for the 
\dfn{thedef}. We closely follow \cite{kac1} and \cite{lz_cqoa}.

\subsubsection{Circle products}
\label{sec:circleprods}
Let $A$ be an algebra. Consider the space of power series 
$A[[z, z\inv]]$. We will write  series $a \in A[[z, z\inv]]$
in the form 
$$
a(z) = \sum_{n\in\Z} a(n) z^{-n-1}, \qquad a(n) \in A.
$$
On $A[[z, z\inv]]$ there is an infinite sequence of bilinear products
$\encirc{n},\ n \in \Z_+$, given by
\begin{equation}\label{fl:circposprod}
\big(a \encirc{n} b\big)(z) = 
\op{Res}_w\big(a(w)\,b(z)\,(z-w)^n\big).
\end{equation}

Explicitly, for a pair of series 
$a(z) = \sum_{n\in\Z} a(n) z^{-n-1}$ and 
$b(z) = \sum_{n\in\Z} b(n) z^{-n-1}$ we have
$$
\big(a \encirc{n} b\big)(z) 
= \sum_m \big(a \encirc{n} b\big)(m) \, z^{-m-1},
$$
where

\begin{equation}\label{fl:explprod}
\big(a \encirc{n} b\big)(m) = 
\sum_{s=0}^n (-1)^s \binom{n}{s}\, a(n-s) b(m+s).
\end{equation}

There is also the linear derivation 
$D = d/dz : A[[z, z\inv]] \to A[[z, z\inv]]$. Easy to see $D$ and 
$\encirc n$ satisfy conditions (ii) and (iii) of \dfn{thedef}.

We can consider formula \fl(explprod) as a system of linear equations
with unknowns $a(k)b(l), \ \ k \in \Z_+, \ l\in \Z$. 
This system is triangular, and its unique solution is given by 
\begin{equation}\label{fl:posproduct}
a(k)b(l) = \sum_{s=0}^k
\binom ks \big(a\encirc s b\big)(k+l-s).
\end{equation}

\begin{Rem}
The term ``circle products'' appears in \cite{lz_cqoa}, where the
product ``$\encirc n$'' is denoted by ``$\circ_n$''.
In \cite{kac2} this product is denoted by ``$\empty_{(n)}$''. 
\end{Rem}

\bigskip
\subsubsection{Locality} 
\label{sec:locality}
Next we define a very important property of power series, which makes them 
form a conformal algebra. Let again $A$ be an algebra.

\begin{Def}
(See \cite{bpz}, \cite{kac_loc}, \cite{kac2}, 
\cite{lz_cqoa}, \cite{lz_preprint})
A series $a \in A[[z, z\inv]]$ is called {\it local of order} $N$ to
$b \in A[[z, z\inv]]$ for some $N \in\Z_+$ if
\begin{equation}\label{fl:locality}
a(w)b(z)\,(z-w)^N = 0.
\end{equation}
If $a$ is local to $b$ and $b$ is local to $a$ then we say that 
$a$ and $b$ are {\it mutually local}.
\end{Def}

\begin{Rem}
In \cite{lz_cqoa} and \cite{lz_preprint} the property \fl(locality) is
called  {\it quantum commutativity}.
\end{Rem}

Note that  \fl(locality) implies that for  every $n\ge N$ one has
$a \encirc n b = 0$.  We will denote 
the order of locality by $N(a,b)$, i.e. 
$$
N(a,b) = \min \{n\in Z_+ \ | \  \forall k \ge n, 
a \encirc k b = 0\}.
$$

Note also that if $A$ is a commutative or skew-commutative algebra, e.g. 
a Lie algebra, then locality is a symmetric relation. In this case we 
say ``$a$ and $b$ are local'' instead of ``mutually local''.

Let $a(z) = \sum_{m\in\Z} a(m) z^{-m-1}$ and 
$b(z) = \sum_{n\in\Z} b(n) z^{-n-1}$ be some series, then 
the locality condition \fl(locality) reads 
\begin{equation}\label{fl:explloc}
\sum_s (-1)^s \binom Ns a(m-s) b(n+s) = 0 \qquad \text{  for any  } n,m \in \Z.
\end{equation}

The locality condition \fl(locality) is known to be equivalent to formula 
\begin{equation}\label{fl:product}
a(m)b(n) = \sum_{s=0}^{N(a,\,b)-1} \binom ms \big(a\encirc s b\big)(m+n-s).
\end{equation}

The following statement is a trivial consequence of the
definitions. 

\begin{Prop}\label{prop:algofser}\sl
Let $A$ be an algebra and let $S \subset A[[z, z\inv]]$ 
be a space of pairwise mutually local
power series, which is
closed under all the circle products and $\partial$. Then $S$ is a 
conformal  algebra. 
\end{Prop}

One can prove (see, for example, \cite{kac2}) that 
such families exhaust all conformal algebras.

Finally, we state here a trivial property of local series:

\begin{Lem}\label{lem:weakdong}\sl
Let $a, b \in A[[z, z\inv]]$ be a pair of formal power series
and assume $a$ is local to $b$. 
Then each of the series $a, Da, 
za$ is local to each of  $b, Db, zb$.
\end{Lem}

\bigskip
\subsection{Construction of the coefficient algebra of a conformal
algebra}
\label{sec:coeff} 
Given a conformal algebra $C$ we can build its {\it coefficient algebra}
$\op{Coeff} C$ in the following way. For each integer $n$ take a linear
space $\^A(n)$ isomorphic to $C$. Let 
$\^A = \bigoplus_{n\in \Z} \^A(n)$. For an element $a\in C$ we will
denote the corresponding element in $\^A(n)$ by $a(n)$. Let $E\subset \^A$
be the subspace spanned by all elements of the form 
\begin{equation}\label{fl:identities}
(Da)(n) + n a(n-1) \qquad
\text{for any  } \  a \in C, \ n \in \Z. 
\end{equation}

The underlying linear space of $\op{Coeff} C$ is $\^A / E $. 
By abuse of notations we will denote the image of $a(n) \in \^A$ in
$\op{Coeff} C$ again by $a(n)$. The following proposition defines the
product on $\op{Coeff} C$. 

\begin{Prop}\label{prop:coeffalg}\sl
Formula \fl(product)
unambiguously defines a bilinear product on $\op{Coeff} C$.
\end{Prop}

Clearly \fl(product) defines a product on $\^A$. To show that the
product is well defined on $\op{Coeff} C$ it is enough to check only
that  
\begin{equation*}
(Da)(m)\, b(n) = -m a(m-1) b(n) \ \text{  and   } \ 
a(m)\, (Db)(n) = -n a(m) b(n-1),
\end{equation*}
which is a straightforward calculation.

\bigskip
\subsection{Examples of conformal algebras}
\label{sec:examples} 
\subsubsection{Differential algebras}\label{sec:diff}
Take a pair $(A, \delta)$, where $A$ is an associative algebra, and 
$\delta:A \to A$ is a locally nilpotent derivation:
$$
\delta(ab) = \delta(a)b + a\delta(b),\qquad \delta^n(a) = 0 \text{ for } n \gg 0.
$$
Consider the ring $A[\delta,\delta\inv]$.  Its elements are polynomials of the
form $\sum_{i\in \Z}a_i \delta^i$,
where only finite number of $a_i \in A$ are nonzero. Here we put 
$a \delta^{-n} = a (\delta\inv)^n$ and $a \delta^0 = a$. 
The multiplication is defined by the formula
$$
a \delta^k\cdot b \delta^l = \sum_{i\ge0}\binom ki a \delta^i(b) \delta^{k+l-i}.
$$
It is easy to check that $A[\delta, \delta\inv]$ is a well-defined associative
algebra. In fact, $A[\delta, \delta\inv]$ is the Ore localization of the ring of
differential operators $A[\delta]$. If in addition $A$ has an identity
element 1, then  $\delta(1)=0$ and $\delta \delta\inv = \delta\inv \delta = 1$. 

For $a\in A$ denote $\~a = 
\sum_{n\in\Z} a \delta^n \, z^{-n-1} \in A[\delta, \delta\inv]\,[[z, z\inv]]$.

One easily checks that 
for any $a,b \in A, \ \~a$ and $\~b$ are local and 
\begin{equation}\label{fl:diffass}
\~a\encirc n \~b = \~{a \delta^n(b)}.
\end{equation}

So by \lem{weakdong} and \prop{algofser} series 
$\{\~a \ | \ a\in A\} \subset A[\delta, \delta\inv]\,[[z, z\inv]]$ 
generate an (associative) conformal algebra, see \sec{varieties}.  

One can instead consider $A[\delta, \delta\inv]$ to be a Lie algebra, with 
respect to the
commutator $\ad pq = pq-qp$. If two series $\~a$ and $\~b$
are local in the associative sense they are  local in the
Lie sense too. One computes also
\begin{equation}\label{fl:difflie}
\~a\encirc n \~b = \~{a \delta^n(b)} - 
\sum_{s\ge0}(-1)^{n+s}\frac 1{s!}\partial^s \(b \delta^{n+s}(a)\)^{\sim}, 
\end{equation}
where $\partial = d/dz$. Note that in \fl(difflie) the circle products
are defined by
\begin{equation}\label{fl:circlie}
\big(\~a \encirc{n} \~b \big)(m) = 
\sum_{s=0}^n (-1)^s \binom{n}{s}\, \ad{a \delta^{n-s}}{b \delta^{m+s}}.
\end{equation}

Again, it follows that 
$\{\~a \ | \ a\in A\} \subset A[\delta, \delta\inv]\,[[z, z\inv]]$ 
generate a (Lie) conformal algebra, see \sec{varieties}.  

An important special case is when there is an element $\upsilon \in A$ such
that $\delta(\upsilon) = 1$. Then 
$\~\upsilon = \sum_n \upsilon \delta^n \, z^{-n-1} \in A[\delta,\delta\inv]\,[[z, z\inv]]$ 
generates with respect to the product \fl(circlie)
a (centerless) Virasoro conformal
algebra. It satisfies the following relations:
$$
\~\upsilon \encirc{0} \~\upsilon = \partial \~\upsilon, \qquad
\~\upsilon \encirc{1} \~\upsilon = 2 \~\upsilon,
$$
and the rest of the products are 0.

\subsubsection{Loop algebras}
\label{sec:loop}
Let $\mathfrak g$ be a Lie algebra over an algebraically closed field 
$\Bbbk$, and let 
$\sigma: \mathfrak g \to \mathfrak g$ be an automorphism of finite
order, $\sigma^p = \op{id}$. Then  $\mathfrak g$ is
decomposed into a direct sum of eigenspaces of $\sigma$:
$$
\mathfrak g = \bigoplus_{k\in \Z/p\Z} \mathfrak g_k,\qquad 
\sigma\big|\raisebox{-5pt}{$\mathfrak g_k$} = e^{2\pi i k/p}.
$$

Define {\it twisted loop algebra} 
$\~{\mathfrak g} \subset \mathfrak g \otimes \Bbbk[t,t\inv]$ by 
$$
\~{\mathfrak g} = 
\Big\{ \sum_j a_j t^j \ \Big| \ 
a_j \in \mathfrak g_{j \op{mod}p}\Big\}.   
$$

The Lie product in  $\~{\mathfrak g}$ is given by 
$\ad{a\otimes t^m}{b\otimes t^n} = \ad ab\otimes t^{m+n}$. If $p=1$,
then $\~{\mathfrak g} = \mathfrak g \otimes \Bbbk[t,t\inv]$, of course.

Now for any $a\in \mathfrak g_k, \ 0\le k < p$, define 
$$
\~a = \sum_{j\in \Z} a t^{pj+k}\, z^{-j-1} \in 
\~{\mathfrak g}\,[[z, z\inv]].
$$

Easy to see that any two $\~a, \~b$ are local with $N(\~a, \~b) = 1$
and if $a\in \mathfrak g_k$ and  $b\in \mathfrak g_l$ we have
$$
\~a\encirc{0} \~b = 
\begin{cases}
\~{\ad ab} & \text{if } k+l < p\\
z\~{\ad ab} & \text{if } k+l \ge p
\end{cases}.
$$

As in \sec{diff}, we conclude that 
$\{\~a \ | \ a\in \mathfrak g\} \subset \~{\mathfrak g}\,[[z, z\inv]]$ 
generate a (Lie) conformal algebra. Again, see \sec{varieties} for
the definition of varieties of conformal algebras. 


\bigskip
\subsection{More on coefficient algebras} 
\label{sec:coeffmore}
Let $C$ be a conformal algebra and let 
$A = \op{Coeff} C$. Define 
\begin{equation*}
\begin{split}
A_+ &= \op{Span}\{ a(n) \ | \ a \in C, \, n \ge 0 \},\\
A_- &= \op{Span}\{ a(n) \ | \ a \in C, \, n < 0 \},\\
A(n) &= \op{Span}\{ a(n) \ | \ a \in C \}.
\end{split}
\end{equation*}
Define also for each $n\in \Z$ linear maps $\phi(n): C \to A(n)$ by 
$a \mapsto a(n)$, and let 
$\phi = \sum_{n\in \Z} \phi(n) z^{-n-1} : C \to A[[z, z\inv]]$ so that 
$\phi a = \sum_{n\in \Z} a(n) z^{-n-1}$.

Here we summarize some  general properties of conformal algebras and their
coefficient algebras.

\begin{Prop}\label{prop:coeff}\sl
\begin{itemize}
\item[(a)]
$A = A_- \oplus A_+$ is a direct sum of subalgebras.
\item[(b)]
$A_+$ and $A_-$ are filtered algebras with filtrations given by
\begin{gather*}
A(0) \subseteq A(1) \subseteq \cdots \subseteq A_+,\qquad
A_- = A(-1) \supseteq A(-2) \supseteq \cdots \\
\bigcup_{n\ge0} A(n) = A_+,\qquad \bigcap_{n<0} A(n) = 0.
\end{gather*}
\item[(c)]
$$
\op{Ker}\phi(n) = 
\begin{cases}
D^{n+1} C + \bigcup_{k\ge1} \op{Ker} D^k & \text{if } \  n \ge 0\\
\op{Ker} D^{-n-1} & \text{if } \ n < 0.
\end{cases}
$$
In particular, $\phi(-1)$ is injective.
\item[(d)]
The map $\phi : C \to A[[z, z\inv]]$, given by 
$a \mapsto \sum_{n\in Z} a(n) z^{-n-1}$ is an injective homomorphism
of conformal algebras, i.e. it preserves the circle products and agrees
with the derivation:
\begin{equation}\label{fl:phi}
\phi(a\encirc n b) = \phi(a) \encirc n \phi(b), \qquad 
\phi(Da) = D \phi(a).
\end{equation}
\item[(e)]
The map  $\phi : C \to A[[z, z\inv]]$ has the following universal
property: for any homomorphism $\psi:C \to B[[z, z\inv]]$ of $C$ to an
algebra of formal power series, there is the unique algebra
homomorphism $\rho:A\to B$ such that the corresponding diagram 
commutes:
\begin{equation*}
\begin{array}{c}
A[[z,z\inv]] \xrightarrow{\ \,\rho\ \,} B[[z, z\inv]] \\
_\phi\!\nwarrow \ \ \ \ \nearrow_\psi \\
 C \\
\end{array}
\end{equation*}
\item[(f)]
The formula $D (a(n)) = -na(n-1)$ defines a derivation $D: A\to A$,
such that $D A_- \subset A_-, \ D A_+ \subset A_+$.
\end{itemize}
\end{Prop}  

\begin{proof}
From formula \fl(product) for the product in $A$ it easily follows that $A_+$
and $A_-$ are indeed subalgebras. Also none of the linear identities
\fl(identities) contain both generators with negative and non-negative
index. This proves (a). Similar arguments establish also (b). 

\medskip

Now we prove that $\op{Ker}\phi(n)$ is included in
the right-hand side of (c). 
Take some $a \in C, \ a \neq 0$, and assume that 
$a(n) = 0$. Then $a(n)\in \^A$ is a linear combination of identities
\fl(identities) (see \sec{coeff}) so we must have in $\^A$ 
$$
a(n) = \sum_{k=k_{\min}}^{k_{\max}} 
\lambda_k \big( (Da_k)(k) + k a_k(k-1)\big).
$$
We can assume that $\lambda_k \neq 0$ for all $k_{\min} \le k \le k_{\max}$ 
and that $a_k \neq 0$ for $k = k_{\min}$
and for $k = k_{\max}$. Assume also that $\lambda_k = 0$ if $k > k_{\max}$ or 
$k < k_{\min}$. 

Comparing terms with index $k$ for $k_{\min} \le k \le k_{\max}$,
we get
\begin{equation}\label{fl:termk}
\delta_{kn} a = \lambda_k Da_k + \lambda_{k+1} (k+1) a_{k+1}. 
\end{equation}

Taking in \fl(termk) $k=k_{\min}-1$ we see that there are two cases: 
either (1)  $k_{\min} =0$ and $n \ge 0$ or (2) $n+1 = k_{\min} \neq 0$.

\medskip\noindent
{\it Case} (1):\quad 
Taking in \fl(termk) $k = 0, \dots, n-1$ we get that 
$a_k \in D^k C$ for $0\le k \le n$. Now we have two subcases:
$k_{\max} > n$ and $k_{\max} \le n$. 

If $k_{\max} > n$ we substitute
in \fl(termk) $k = k_{\max},  k_{\max}-1, \dots, n+1$ and get that
$D^{k_{\max} -k + 1} a_k = 0$. Now take $k=n$ in  \fl(termk) and get 
that $a \in D^{n+1} C + \op{Ker}D^{k_{\max}-n}$.

If $k_{\max} \le n$ we have $\lambda_{n+1} = 0$, and hence substitution
$k=n$ in \fl(termk) gives $a \in D^{n+1} C$.

\medskip\noindent
{\it Case} (2):\quad 
Here we again have two subcases: $n \ge 0$ or $n < 0$.

If $n \ge 0$ then as in the previous case, 
we get $D^{k_{\max} -k + 1} a_k = 0$ for $n+1\le k\le k_{\max}$.
Now taking $k=n$ in \fl(termk) we get $a \in \op{Ker}D^{k_{\max}-n}$.

Finally, if $n < 0$ then, since $\lambda_n = 0$, we have
$a = \lambda_{n+1}(n+1) a_{n+1}$. Then 
we substitute $k = n+1, n+2, \dots \ $ into 
\fl(termk) until for some $k_0\le -1$ we get 
$\lambda_{k_0+1}(k_0+1) a_{k_0+1} = 0$.  It  follows that 
$D^{k_0-k+1}a_k = 0$ for $n+1 \le k \le k_0$, therefore
$a \in \op{Ker} D^{k_0-n} \subset \op{Ker} D^{-n-1}$. 
This proves one inclusion in (c).
It also follows that $\op{Ker}\phi(-1) = 0$.

\medskip

Next we show that $\phi$ is a homomorphism of conformal algebras, that
is, formulas \fl(phi) hold. For the first 
identity we have
$$
\phi(Da) = \sum_{n\in\Z} (Da)(n) z^{-n-1} = 
\sum_{n\in\Z} (-n) a(n-1) z^{-n-1} = 
\frac d{dz} \sum_n a(n) z^{-n-1}.
$$

The second identity reads 
\begin{equation*}
\big(a\encirc n b\big)(m) = \sum_s (-1)^s \binom ns a(n-s) b(m+s),
\end{equation*}
which is precisely the formula \fl(explprod).

Now (d) is done after we notice that $\phi$ is injective, since 
$\phi(-1)$ is injective.

\medskip
Now we can prove the other inclusion in (c). If $a \in \op{Ker}D^k C$, 
then $\phi a$ is a solution of differential equation 
$\partial_z^k \phi a(z) =0$, hence $\phi a$ is a polynomial of degree
at most $k-1$, therefore $\phi(n) a = 0$ for $n \ge 0$ and for 
$n < -k$. If $a \in D^k C$, then $\phi(n) a = 0$ for $0\le n \le k-1$, 
by induction and \fl(identities).


\medskip
Statement (e) is clear, since identities \fl(identities) hold for any
homomorphism $\psi: C \to B[[z, z\inv]]$.

\medskip
Finally, the formula $D(a(n)) = -n a(n-1)$ defines a derivative of
$\^A$. So in order to prove (f) we have to show that $D$ agrees with
the identities \fl(identities). This is indeed the case:
$$
D\big((Da)(n)+ n a(n-1)\big) = -n\big((Da)(n-1) + (n-1) a(n-2)\big).
$$
\end{proof}

\bigskip
\subsection{Varieties of conformal  algebras}
\label{sec:varieties} 
Consider now a variety of algebras $\mathfrak A$ (see \cite{cohn}, 
\cite{jacobson}). 
\begin{Def}\label{dfn:confA}
A conformal algebra $C$ is a {\it $\mathfrak A$-conformal algebra}
if  $\op{Coeff} C$ lies in the variety $\mathfrak A$.
\end{Def}

The identities in $\mathfrak A$-conformal algebras are
all the circle-products identities $R$ such that for any integer $m$, \ 
$R(m)$ becomes an  $\mathfrak A$-algebra identity after substitution of
\fl(explprod) for every circle product in $R$. Conversely, 
given a classical algebra identity $r$ we can substitute \fl(product)
for all  products in $r$ and get an 
identity of $\mathfrak A$-conformal algebras. This way we get a
correspondence between classical and conformal identities.
See the next section for examples.

Combining \prop{algofser} and (d) of \prop{coeff} we get
the following well-known fact:
\begin{Prop}\label{prop:confrealis}\sl
$\mathfrak A$-conformal algebras  are exhausted (up to isomorphism) 
by conformal 
algebras of formal power series $S \subset A[[z,z\inv]]$ for 
$\mathfrak A$-algebras $A$.
\end{Prop}

\bigskip
\subsection{Associative and Lie conformal  algebras}
\label{sec:asslie} 
The following theorem gives the explicit correspondence between 
conformal and classical algebras in some important cases.

\begin{Thm}[See \cite{kac_fd}]\label{thm:identity}\sl
Let $C$ be a conformal algebra and $A = \op{Coeff} C$ its coefficient
algebra. 
\begin{itemize}
\item[(a)]
$A$ is associative if and only if the following identity holds in $C$:
\begin{equation}\label{fl:assconf}
(a \encirc{n} b)\encirc m c = 
\sum_{s=0}^n (-1)^s\binom{n}{s} 
a\encirc{n-s} (b \encirc{m+s} c).
\end{equation}
\item[(b)]
The Jacoby identity 
$\ad{\ad{a}{b}}{c} = \ad{a}{\ad{b}{c}} - \ad{b}{\ad{a}{c}}$ 
in $A$ is equivalent to the following conformal Jacoby identity in $C$:
\begin{equation}\label{fl:jacconf}
(a \encirc{n} b)\encirc m c = 
\sum_{s=0}^n (-1)^s\binom{n}{s} 
\big(a\encirc{n-s} (b \encirc{m+s} c) -
b\encirc{m+s} (a\encirc{n-s} c)\big)
\end{equation}
\item[(c)]
The skew-commutativity  identity 
$\ad{a}{b} = -\ad{b}{a}$ in $A$
corresponds to the quasisymmetry identity:
\begin{equation}\label{fl:quasisym}
a\encirc n b = \sum_{s\ge0} (-1)^{n+s+1}\frac 1{s!} D^s (b\encirc{n+s}a).
\end{equation}
\item[(d)]
The commutativity of $A$ is equivalent to 
\begin{equation}\label{fl:comconf}
a\encirc n b = \sum_{s\ge0} (-1)^{n+s}\frac 1{s!} D^s (b\encirc{n+s}a)
\end{equation}
\end{itemize}
\end{Thm}

The identities \fl(assconf), \fl(jacconf) and \fl(quasisym) immediately
imply the following
\begin{Cor}\label{cor:module}\sl
Let $C$ be a Lie conformal or an associative conformal algebra, 
and $A = \op{Coeff} C$ its coefficient algebra. Then
$C$ is an $A_+$-module with the action given by $a(n)c = a\encirc n \!c\,$
for $a,c \in C, \, n  \in \Z_+$.
Moreover, this action agrees with the derivations on $A_+$ and $C$: \  
$(D a(n)) c = \ad D{a(n)}c$.
\end{Cor}

From now on we will deal only with associative or Lie conformal 
algebras.

\bigskip
\subsection{Dong's lemma}
\label{sec:dong}
We end this section by stating a very important property of formal
power series over 
associative or Lie algebras. This property  allows to construct
conformal algebras by taking a collection of generating series.

\begin{Lem}\label{lem:dong}\sl
Let $A$ be an associative or a Lie algebra, and 
let $a,\,b,\,c \in A[[z,z\inv]]$ be three formal power series.
Assume that they are pairwise mutually local. 
Then for all $n \in \Z_+,\ a \encirc n b$ and $c$ are mutually
local. Moreover, in the Lie algebra case, 
\begin{equation}\label{fl:estimate}
N(a\encirc n b, c) = N(c, a\encirc n b) \le N(a,b) + N(b,c) + N(c,a)-n-1,
\end{equation}
and in the associative case
$$
N(a\encirc n b, c) \le N(b,c),\qquad 
N(c, a\encirc n b) \le N(c,a)+N(a,b)-n-1.
$$
\end{Lem}

\bigskip

\section{Vertex algebras}
\subsection{Fields}\label{sec:fields}
Let now $V$ be a vector space over $\Bbbk$. Denote by $gl(V)$ the Lie
algebra of all $\Bbbk$-linear operators on $V$.
Consider the space 
$\F(V) \subset gl(V) [[z, z^{-1}]]$
of {\it fields} on $V$, given by
$$
\F(V) = \left\{ \left.
\sum_{n\in \mathbb Z}a(n)\,z^{-n-1}\ \right|
\forall v \in V, \  a(n) v = 0 \text{  for } n \gg 0
\right\}.
$$
For $a(z) \in \F(V)$ denote 
$$
a_-(z) = \sum_{n<0}a(n)\,z^{-n-1},\quad
a_+(z) = \sum_{n\ge 0}a(n)\,z^{-n-1}.
$$
Denote also  by  $ \1 = \1_{\F(V)} \in \F(V)$ the identity operator, 
such that $ \1 (-1) = \op{Id}_V$,  all other coefficients are 0. 

\begin{Rem}
In \cite{lz_cqoa} and \cite{lz_preprint} the elements of 
$\F(V)$ are called {\it quantum operators} on $V$.
\end{Rem}

We view $gl(V)$ as a Lie algebra, and \sec{circleprods} gives 
a collection of products $\encirc n, \ n \in \Z_+,$ on $\F(V)$. Now
in addition to these products 
we introduce products $\encirc n$  for  $n<0$. Define first 
$\encirc{-1}$ by 
\begin{equation}\label{fl:minus1}
a(z) \encirc{-1} b(z) =  a_-(z) b(z) + b(z) a_+(z).
\end{equation}
Note that the products in \fl(minus1) make sense, since for any 
$v \in V$ we have $a(n)v = b(n)v = 0$ for $n \gg 0$.
The  $-1$-st product is also known as 
the {\it normally ordered product} (or {\it Wick  product}) and
is usually denoted by $\:a(z)b(z)\:$. 

Next, for any $n<0$ set 
\begin{equation}\label{fl:negprod}
a(z) \encirc{n} b(z) = 
\tfrac1{(-n-1)!}\, \: \big(D^{-n-1}a(z)\big) b(z) \:,
\end{equation}
where $D = \frac d{dz}$. 
Taking $b = \1$ we get 
\begin{equation}\label{fl:minus2}
a \encirc{-1} \1 = a, \qquad a \encirc{-2} \1 = D a. 
\end{equation}
Easy to see that 
$$
\1 \encirc n a = \delta_{-1,n} a.
$$

We have the following explicit formula for the circle products:
if $\big(a \encirc{n} b\big)(z)
 = \sum_m \big(a \encirc{n} b\big)(m) z^{-m-1}$, then 
\begin{equation}\label{fl:expl}
\begin{split}
\big(a \encirc{n} b\big)(m) 
& = \sum_{s\le n} (-1)^{s+n} \binom n{n-s} a(s)b(m+n-s)\\
& - \sum_{s\ge0} (-1)^{s+n} \binom n{s} b(m+n-s)a(s).
\end{split}
\end{equation}

Note that if $n>0$ then \fl(expl) becomes
$$
\big(a \encirc{n} b\big)(m) = \sum_{s\ge0} (-1)^{n+s} \binom ns
\ad{a(s)}{b(m+n-s)},
$$
which is precisely formula \fl(explprod) for Lie algebras.

It is easy to see, that $D$ is a
derivation of all the circle products:
\begin{equation}\label{fl:der}
D(a\encirc n b) = D a \encirc n b + a \encirc n D b.
\end{equation}

Note also that the Dong' \lem{dong} remains valid for negative $n$
and the estimate \fl(estimate) still holds.

\bigskip
\subsection{Definition of vertex algebras}
\label{sec:vertex}
Instead of giving a formal definition  of
vertex algebra in spirit of \dfn{thedef},
we present a description of these algebras
similar to \prop{confrealis}. For a more abstract approach see e.g.
\cite{fhl}, \cite{kac2}, \cite{lz_cqoa} or \cite{lz_preprint}.

\begin{Def}\label{vertex}
A {\it vertex algebra} is a subspace $S\subset \F(V)$ of fields over a
vector space $V$ such that 
\begin{itemize}
\item[(i)]
Any two fields $a,b \in S$ are local (in the Lie sense).
\item[(ii)]
$S$ is closed under all the circle products $\encirc n, \ n \in \Z$, given by 
\fl(expl).
\item[(iii)]
$\1 \in S$.
\end{itemize}
\end{Def}
Note that from \fl(minus2) it follows that a vertex algebra is closed under
the derivation $D = d/dz$.

Note also that a vertex algebra is a Lie conformal algebra. 

Let $S \subset \F(V)$ be a vertex algebra. We introduce the left
action map $Y: S \to \F(S)$ defined by  
\begin{equation}\label{fl:Y}
Y(a) = \sum_{n\in \Z} (a \encirc{n} \cdot\ )\, \zeta^{-n-1}. 
\end{equation}
Clearly, $Y(\1_S) = \1_{\F(S)}$. 

We state here the following characterizing property of $Y$ 
(see \cite{kac2} or \cite{lz_preprint}):
\begin{Prop}\label{prop:vertex}\sl
The left action map $Y: S \to \F(S)$ is an isomorphism of vertex
algebras, i.e. $Y(S) \subset \F(S)$ is a vertex algebra and 
\begin{equation}\label{fl:Y_is_hom}
Y(a\encirc n b) = Y(a)\encirc n Y(b), \qquad
Y(\1_S) = \1_{\F(S)}. 
\end{equation}
\end{Prop}

From \fl(minus2) and \fl(der) it follows that $Y$ also 
agrees with $D$: 
$$
Y(Da) = \partial_\zeta Y(a) = \ad D{Y(a)}.
$$

\bigskip
\subsection{Enveloping vertex algebras of a Lie conformal algebra}
\label{sec:envertex}
Let  $C$ be a Lie conformal algebra and $L=\op{Coeff}C$ its 
coefficient Lie algebra.

\begin{Def}(See \cite{kac1}, \cite{kac2})
\begin{itemize}
\item[(a)]
An $L$-module $M$ is called {\it restricted} if for any $a\in C$ and 
$v \in M$ there is some integer $N$ such that for any $n\ge N$ one has 
$a(n) v = 0$.
\item[(b)]
An $L$-module $M$ is called a {\it highest weight} module if 
it is generated over $L$ by a single element $m \in M$ such that 
$L_+ m = 0$. In this case $m$ is called the {\it highest weight vector}
\end{itemize}
\end{Def}
Clearly any submodule and any factor-module of a restricted module are
restricted.

Let $M$ be a restricted $L$-module. Then the representation
$\rho: L \to gl(M)$ could be extended to the map 
$\rho: L[[z, z\inv]] \to \F(M)$ which combined with the canonical
embedding $\phi:C \to L[[z, z\inv]]$ (see (d) of \prop{coeff}) gives
conformal algebra homomorphism $\psi: C \to \F(M)$. 
Then $\psi(C)\subset \F(M)$ consists of pairwise local fields, and by 
Dong's \lem{dong},  $\psi(C)$ together with $\1 \in \F(M)$ generates
a vertex algebra $S_M\subset \F(M)$. 

The following proposition is well-known, see e.g. \cite{frzh}.
\begin{Prop}\label{prop:module}\sl
\begin{itemize}
\item[(a)]
The vertex algebra $S = S_M$ has a structure of a highest weight module 
over $L$ with the highest weight vector $\1$. The action is given by
$$
a(n)\beta = \psi(a) \encirc n \beta, \qquad 
a\in C, \ n \in \Z, \ \beta \in S_M.
$$
Moreover this action agrees with the derivations: 
$$
(Da(n))\beta = \ad D{a(n)} \beta.
$$
\item[(b)]
Any $L$-submodule of $S$ is a vertex algebra ideal. If $M_1$ and $M_2$
are two restricted $L$-modules, $S_1 = S_{M_1}, \ S_2 = S_{M_2}$, and
$\mu:S_1 \to S_2$ is an $L$-module homomorphism such that $\mu(\1) = \1$, 
then $\mu$ is a vertex algebra homomorphism. 
\end{itemize}
\end{Prop}

\bigskip
\subsection{Universal enveloping vertex algebras}
\label{sec:univertex}
Now we build a universal highest weight module $V$ over $L$, which is
often referred to as {\it Verma module}.  Take the 1-dimensional 
trivial $L_+$-module $\Bbbk \1_V$, generated by an element $\1_V$. 
Then let 
$$
V = \op{Ind}_{L_+}^L \Bbbk \1_V = U(L) \otimes_{U(L_+)}\Bbbk \1_V
\cong U(L)/U(L)L_+.
$$ 
It is easy to see that $V$
is a restricted module and hence we get an enveloping vertex algebra 
$S=S_V\subset \F(V)$ and an homomorphism $\psi:C\to S$. Clearly, 
$\psi$ is injective, since $\rho:L\to gl(V)$ is injective. 

\begin{Thm}\label{thm:univertex}\sl
\begin{itemize}
\item[(a)]
The map $\chi:S \to V$ given by $\alpha \mapsto \alpha(-1)\1_V$
is an $L$-module isomorphism, and $\chi(\1_{S}) = \1_V$.
\item[(b)]
$S$ is the universal enveloping vertex algebra of $C$ in the following
sense:
If $\mu:C \to U$ is another homomorphism of $C$ to a vertex algebra 
$U$, then there is the unique map $\^\mu:S\to U$ which makes up the  
following commutative triangle:
\begin{equation*}
\begin{array}{c}
S \xrightarrow{\ \ \^\mu\ \ } U \\
_\psi\!\nwarrow \ \ \ \nearrow_\mu \\
 C \\
\end{array}
\end{equation*}
\end{itemize}
\end{Thm}

From now on we identify $V$ and $S=S_V$ via $\chi$ and write $V = V(C)$
for the universal enveloping vertex algebra of a Lie conformal algebra 
$C$ and $\1_S = \1_V = \1$. 
The embedding $\psi: C \to V = U(L)/U(L)L_+$ is then given by 
$a\mapsto a(-1)\1$. 
By (c) of \prop{coeff}, the map $\phi(-1):C \to L_-$, defined by
$a\mapsto a(-1)$, is an isomorphism of linear spaces. Therefore,
the image of $C$ in $V$ is equal to $\psi(C) = L_-\1 = L\1 \subset V$.

\bigskip

\section{Free conformal algebras}
\subsection{Definition of free conformal and free vertex algebras}
\label{sec:freeconf}
Let $\mathcal B$ be a set of symbols. Consider a 
function $N: \mathcal B \times  \mathcal B \to \Z_+$,
which will be called a {\it locality function}.

Let $\mathfrak A$ be a variety of algebras. In all the applications 
$\mathfrak A$ will be either Lie or associative algebras.
Consider the category $\mathfrak{Conf}(N)$ of  
$\mathfrak A$-conformal algebras (see \sec{varieties})
generated by the set $\mathcal B$ such that in any conformal algebra
$C \in {\mathfrak Conf}(N)$ one has
$$
a\encirc n b = 0 \qquad
\forall a,b \in \mathcal B \ \ \forall n \ge N(a,b).
$$
By abuse of notations we will not make a distinction between $\mathcal B$
and its image in a conformal algebra $C \in  {\mathfrak Conf}(N)$. 

The morphisms of $\mathfrak{Conf}(N)$ are, naturally, conformal
algebra homomorphisms $f:C \to C'$ such that $f(a) = a$ for any 
$a \in \mathcal B$. 

We claim that $\mathfrak{Conf}(N)$ has the universal object,
a conformal algebra $C=C(N)$, such that for any other 
$C' \in \mathfrak{Conf}(N)$ there is the unique morphism $f:C \to C'$.
We call $C(N)$ a {\it free conformal algebra}, corresponding to the
locality function $N$.

In order to build $C(N)$, we first build the corresponding coefficient
algebra $A = \op{Coeff}C$ (see \sec{coeff}).

Let $A \in \mathfrak A$ be the algebra presented by the set of
generators
\begin{equation}\label{fl:gen}
X = \{b(n)\ | \ b\in \mathcal B, \ n \in \Z \}
\end{equation}
with relations 
\begin{equation}\label{fl:rel}
\left\{\left.
\sum_s (-1)^s \binom {N(b,a)}s b(n-s) a(m+s) = 0 \ \right| \ 
a,b \in \mathcal B, \ \ m,n \in \Z \right\}.
\end{equation}
For any $b \in \mathcal B$ let $\~b = \sum_n b(n) z^{-n-1} \in  A[[z,z\inv]]$.
From \fl(rel) follows that any two $\~a$ and  $\~b$ are mutually local, 
therefore by Dong's \lem{dong} they generate a conformal algebra
$C \subset A[[z,z\inv]]$.  

\begin{Prop}\label{prop:freeconf}\sl
\begin{itemize}
\item[(a)]
$A = \op{Coeff}C$.
\item[(b)]
The conformal algebra $C$ is the free conformal algebra 
corresponding to the locality function $N$.
\end{itemize}
\end{Prop}

\begin{proof}
(a)
Clearly, there is a surjective homomorphism $A \to \op{Coeff}C$, since
relations \fl(rel) must hold in $\op{Coeff}C$. Now the claim follows from
the universal property of $\op{Coeff}C$ (see (e) of \prop{coeff}). 

\medskip\noindent
(b)
Take another algebra $C' \in \mathfrak{Conf}(N)$, and let 
$A' = \op{Coeff}C'$.
Obviously, there is an algebra homomorphism $f:A \to A'$ such that 
$f(b(n)) = b(n)$ for any $b\in \mathcal B$ and $n \in \Z$. It could be 
extended to a map $f: A[[z, z\inv]] \to A'[[z, z\inv]]$. Now it is
easy to see that the restriction $f|\raisebox{-5pt}{$C$}$  gives the 
desired conformal algebra homomorphism $C \to C'$:
$$
\begin{CD}
A[[z,z\inv]] @>f>> A'[[z,z\inv]] \\
@AA A           @AA A \\
C @> f >>  C'
\end{CD}
$$
Indeed, due to
formula \fl(explprod), $f$ preserves the circle products, and, since
$\partial$ is a derivation  of the products, and $f(\partial \~a) =
\partial f(\~a)$,  for $a\in C$ one has  $f(\partial \phi) = \partial
f(\phi)$ for any $\phi \in C$. 
\end{proof}

In case when $\mathfrak A$ is the variety of Lie algebras, 
we may consider the universal vertex enveloping
algebra $V(C)$ of a free Lie conformal algebra $C = C(N)$.
In accordance with \thm{univertex}, we call $V(C)$  a 
{\it free vertex algebra}.

Though the construction of a free conformal and vertex algebras
makes sense for an arbitrary locality function 
$N : \mathcal B \times \mathcal B \to \Z_+$, results of
\sec{basisver}--\sec{basislie}
are valid only for the case when $N$ is  constant.

\bigskip
\subsection{The positive subalgebra of $\op{Coeff} C(N)$}
\label{sec:positive}
Let again $C = C(N)$ be a free conformal algebra corresponding to a
locality function $N: \mathcal B \times \mathcal B \to \Z_+$, 
$\mathcal B$ being an alphabet, and let $A = \op{Coeff}C$.
Recall that by \prop{coeff} (a)
we have the decomposition $A = A_-\oplus A_+$ of the coefficient algebra
into the direct sum of two subalgebras. Denote 
$X_i = \{ b(n) \ | \ b\in \mathcal B, n \ge i \} \subset X$.

\begin{Lem}\label{lem:Aplus}\sl
The subalgebra $A_+\subset A$ is isomorphic to the algebra
$\^A_+$ presented by the set of generators $X_0$ and
those of relations \fl(rel) which contain only elements of
$X_0$:
\begin{equation}\label{fl:rel0}
\left\{\left.
\sum_s (-1)^s \binom {N(b,a)}s b(n-s) a(m+s) = 0 \ \right| \ 
a,b \in \mathcal B, \  m\ge 0, \ n \ge N(b,a) \right\}.
\end{equation}
\end{Lem}
\begin{proof}
Clearly, there is a surjective homomorphism $\varphi:\^A_+ \to A_+$
which maps $X_0$ to itself. We prove that $\varphi$ is in fact 
an isomorphism. We proceed in four steps.

\medskip\noindent{\it Step 1.}\hskip10pt
First we prove that $A_+$ is generated by $X_0$ in $A$. Indeed, 
we have $X_0 \subset A_+$. 
On the other hand, $A_+$ is spanned by elements of the form $a(m)$, where
$m\ge 0$ and $a \in C$ is a circle-product monomial in $\mathcal B$.
By induction on the length of $a$ it is enough to check that if
$a = a_1 \encirc k a_2$, then $a(m)$ is in the subalgebra, generated
by $X_0$, which follows from \fl(explprod).

\medskip\noindent{\it Step 2.}\hskip10pt
Let $\^\tau:\^A_+ \to \^A_+$ be the homomorphism, which acts
on the generators $X_0$ by $a(n) \mapsto a(n+1)$, 
so that $\^\tau(\^A_+)$ is the subalgebra
of $\^A_+$ generated
by $X_1$. We claim that $\^\tau$ is injective, and therefore
$\^\tau(\^A_+) \cong \^A_+$. Indeed, $\^\tau$ acts on the free associative
algebra $\Bbbk\<X_0\>$. Assume that for some 
$p \in \^A_+$ we have $\^\tau(p) = 0$. Take any preimage
$P \in \Bbbk\<X_0\>$ of $p$. Then we have  
$\^\tau(P) = \sum_i \xi_i R_i$, where $\xi_i \in \Bbbk\<X_0\>$
and $R_i$ are relations \fl(rel0), such that in all $\xi_i$ and 
$R_i$ appear only indexes greater or equal to 1. But then
$P$ itself must be of the form $\sum_i \xi'_i R'_i$, where $"\,'\ ``$ 
stands for decreasing all indexes by 1, hence $p = 0$.

\medskip\noindent{\it Step 3.}\hskip10pt
Next we claim that there is an automorphism $\tau$ of the algebra
$A$ which acts on the generators $X$ by the shift
$a(n) \mapsto a(n+1)$. Indeed, relations \fl(rel) are invariant 
under the shift, and clearly, $\tau$ is invertible. 
For any integer $n$  denote $A_n = \tau^n A_+$. We have 
$A_n \cong A_+ = A_0$ for every $n$.

\medskip\noindent{\it Step 4.}\hskip10pt
Now for each integer $n$ take a copy $\^A_n$ of $\^A_+$. Let 
$\^\tau_n : \^A_n \to \^A_{n-1}$ be the isomorphism 
of $\^A_+$ onto $\^\tau(\^A_+)$, built in Step 1. 
Let $\^A$ be the limit of all
these $\^A_n$ with respect to the maps $\^\tau_n$. 
We identify generators of $\^A_n$ with the set $X_n$. 
It is easy to see that $\varphi: \^A_0 \to A_0$ extends to the homomorphism
$\varphi: \^A \to A$, such that $\varphi(\^A_n) = A_n$ and
$\varphi\big|\raisebox{-5pt}{$X$} = \op{id}$. 
Now we observe that all the defining
relations \fl(rel) of $A$ hold in $\^A$, hence there is an 
inverse map $\varphi\inv : A \to \^A$, and therefore
$\varphi$ is an isomorphism.
\end{proof}

\bigskip
\subsection{The Diamond Lemma}
\label{sec:diamond}
For the future purposes
we need a digression on the Diamond Lemma for associative algebras. 
We closely follow \cite{bergman}, but use more modern terminology.  

Let $X$ be some alphabet and $K$ be some commutative
ring. Consider the free associative algebra $K\<X\>$
of non-commutative polynomials with coefficients in $K$.
Denote by $X^*$ the set of words in $X$, i.e. the free
semigroup with 1 generated by $X$.

A {\it rule} on $K\<X\>$ is a pair 
$\rho = (w, f)$, consisting of a word
$w \in X^*$  and a polynomial
$f \in K\<X\>$.
The left-hand side $w$ 
is called {\it the principal part} of rule $\rho$. We will denote 
$w = \bar \rho$.

Let $\mathcal R$ be a collection of rules on $K\<X\>$.
For a rule $\rho = (w,f)\in \mathcal R$ 
and  a pair of words $u,v \in X^*$
consider the $K$-linear endomorphism 
$r_{u\rho v}: K\<X\> \to  K\<X\>$, which fixes
all words in $X^*$ except for $u w v$, and sends the latter
to $u f v$. 

A rule $\rho = (w, f)$ is said to
be {\it applicable} to a word $v \in X^*$ if 
$w$ is a subword of $v$, i.e. $v = v' w v''$. 
The result of application of $\rho$ to $v$ is, naturally,
$r_{v'\rho v''}(v) = v'fv''$.
If $p \in K\<X\>$ is a polynomial which involves a word
$v$, such that a rule $\rho$ is applicable to $v$, then 
we say that $\rho$ is applicable to $p$. 

A polynomial $p\in K\<X\>$ is called {\it terminal} if no rule from 
$\mathcal R$ is applicable to $v$, that is, 
no term of $p$ is of the form $u\bar\rho v$ for $\rho \in \mathcal R$.

Define a binary relation "$\ra$" on $K\<X\>$ in the following way:
Set $p\ra q$ if and only if there is a finite sequence of rules 
$\rho_1, \dots, \rho_n \in \mathcal R$, 
and a pair of  sequences of words $u_i, v_i \in X^*$ such that 
$q = r_{u_n\rho_n v_n} \cdots  r_{u_1\rho_1 v_1} (p)$.

\begin{Def}\label{dfn:rewriting}
\begin{itemize}
\item[(a)]
A set or rules $\mathcal R$ is a {\it rewriting system } on 
$K\<X\>$ if  there are no infinite sequences of
the form 
$$
p_1 \ra p_2 \ra \dots,
$$
i.e. any polynomial $p \in K\<X\>$
can be modified only finitely many times by rules from $\mathcal R$.
\item[(b)]
A rewriting system is {\it confluent} if for any polynomial 
$p \in K\<X\>$ there is the unique terminal polynomial
$t$ such that $p \ra t$.
\end{itemize}
\end{Def}

Any rule $\rho = (w, f) \in \mathcal R$ gives rise to an 
identity $w - f \in  K\<X\>$. Let 
$I(\mathcal R) \subset K\<X\>$
is the two-sided ideal generated by all such  identities.

Let $v_1,v_2 \in X^*$ be a pair of words. A word 
$w \in X^*$ is called {\it composition} of $v_1$ and $v_2$
if $w = w'uw'', \ v_1 = w'u, \ v_2 = uw''$ and $u \neq 0$. 

Finally, take a word $v \in X^*$.
Let us call it an {ambiguity} if there
are two rules $\rho, \sigma \in \mathcal R$ such that either $v$ is a
composition of $\bar\rho$ and $\bar\sigma$ or if $v = \bar\rho$ and
$\bar\sigma$ is a subword of $\bar\rho$. 

Now we can state the Lemma.
\begin{Lem}[Diamond Lemma]\label{lem:diamond}\sl
\begin{itemize}
\item[(a)]
A rewriting system $\mathcal R$ is confluent if and only if
all terminal monomials form a basis of 
$K\<X\>/I(\mathcal R)$.
\item[(b)]
A rewriting system is confluent if and only if
it is confluent on all the ambiguities, that is, 
for any ambiguity $v\in X^*$ there is the unique terminal 
$t\in K\<X\>$ such that  $v \ra t$.
\end{itemize}
\end{Lem}

\begin{Rem}
Statement (a) appears in \cite{newman}. A variant of \lem{diamond}
appears in \cite{bokut72} and \cite{bokut76}. It was also known to
Shirshov (see \cite{shirshov62b}). The name ``Diamond'' is 
due to the following graphical description of the confluency property,
see \cite{newman}.
Let $\mathcal R$ be a rewriting system in sense of \dfn{rewriting} (a),
and let ``$\ra$'' be defined as above. Assume 
 $p,q_1, q_2 \in K\<X\>$ are such that  $p\ra q_1$ and $p\ra q_2$. Then 
there is some $t \in  K\<X\>$ such that $q_1 \ra t$ and $q_2 \ra t$:
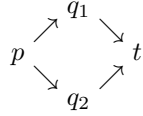
\begin{equation*}
\begin{array}{@{\extracolsep{-7pt}}ccccc}
 &         &q_1&        & \\[-3pt]
 &\nearrow &   &\searrow& \\[-3pt]
p&         &   &        &t\\[-3pt]
 &\searrow &   &\nearrow& \\[-4pt]
 &         &q_2&        & 
\end{array}
\end{equation*}
\end{Rem}

G. Bergman in \cite{bergman} uses the existence of a semigroup order
with descending chain condition on the set of words $X^*$. Though in 
our case there is an order on the set \fl(gen), this order does not
satisfy the descending chain condition, so we slightly modify 
the argument in \cite{bergman}.

\begin{proof}[Proof of \lem{diamond}]
(a)
Assume that the rewriting system $\mathcal R$ is confluent. Define a map
$r:K\<X\> \to K\<X\>$ by taking $r(p)$ to be the unique terminal 
monomial such that $p \ra r(p)$.
The crucial observation is that $r$ is a $K$-linear endomorphism of
$K\<X\>$. So if 
$p = \sum_i \xi_i u_i (w_i - f_i) v_i \in I(\mathcal R),\ \xi_i \in K, \ 
u_i, v_i \in X^*, \ (w_i, f_i) \in \mathcal R$, 
then $r(p) = \sum_i \xi_i r(u_i (w_i - f_i) v_i) = 0$, therefore 
the terminal monomials are linearly independent modulo $I(\mathcal R)$.

From the other side, if $\mathcal R$ is not confluent, then there are
a polynomial $p \in K\<X\>$ and terminals $q_1, q_2 \in K\<X\>$  
such that $p \ra q_1, \ p\ra q_2$
and $q_1 \neq q_2$, and then $q_1 - q_2 \in I(\mathcal R)$.

\medskip\noindent
(b)
Take a polynomial $p\in K\<X\>$. We prove that there is the unique
terminal $t$ such that $p\ra t$ by induction on the number 
$n(p) = \#\{q\ | \ p\ra q\}$. The condition (a) of \dfn{rewriting}
assures that $n(p)$ is always finite.

If $n(p) = 0$ then $p$ is a terminal itself and there is nothing to 
prove. By induction, without loss of generality we can assume
that there are at least two different rules $\rho, \sigma \in \mathcal R$
which are applicable to $p$. It means that there are some
words $u,v,x,y \in X^*$ such that $r_{u\rho v} (p) \neq p$, 
$r_{x\sigma y} (p) \neq p$ and 
$r_{u\rho v} (p) \neq r_{x\sigma y} (p)$.
By induction, both $r_{u\rho v} (p)$ and $r_{x\sigma y} (p)$ are uniquely
reduced to terminals, say 
$r_{u\rho v} (p)\ra t_1$ and $r_{x\sigma y} (p)\ra t_2$. We need to show 
that $t_1 = t_2$.

Consider two cases: 
when $\bar \rho$ and $\bar\sigma$ have common symbols in $p$, and 
thus $u\bar\rho v = x \bar\sigma y$ is  a word in $p$, and when 
$\bar \rho$ and $\bar\sigma$ are disjoint. 

In the first case, let $w \in X^*$ be the union of $\bar \rho$ and 
$\bar\sigma$ in $p$. Then $w$ is an ambiguity. By assumption, 
there is the unique terminal $s \in K\<X\>$ such that $w \ra s$.
Let $q  \in K\<X\>$ be obtained from $p$ by substituting $w$ by $s$.
Then we have
\begin{equation}\label{fl:diagram}
\begin{array}{@{\extracolsep{-7pt}}ccccc}
 &         &r_{u\bar\rho v}(p)  &        & \\
 &\nearrow &                    &\searrow& \\
p&         &                    &        &q\\
 &\searrow &                    &\nearrow& \\
 &         &r_{x\bar\sigma y}(p)&        & \\
\end{array}
\end{equation}
By induction, $q$ is uniquely reduced to a terminal $t$, therefore one
has $r_{u\rho v} (p)\ra t$ and $r_{x\sigma y} (p)\ra t$.

In the second case note that 
$r_{x\sigma y}r_{u\rho v}(p) = r_{u\rho v}r_{x\sigma y}(p)$. Denote
this polynomial by $q$. Then relations \fl(diagram) still hold, and we
finish by the same argument as in the first case. 
\end{proof}

\bigskip
\subsection{Basis of a free vertex algebra}
\label{sec:basisver} 
Return to the setup of \sec{freeconf}.
From now on we take the locality function  $N(a,b)$ to be
constant: $N(a,b) \equiv N$. Let $C = C(N)$ be the free Lie
conformal algebra and $L = \op{Coeff} C$ its  Lie
algebra of coefficients, see \prop{freeconf}. In this section we build
a basis of the universal enveloping algebra $U(L)$ of $L$ 
and a basis of the free vertex algebra $V = V(C)$.

We start with endowing  $\mathcal B$  with an arbitrary linear order.   
Then define a linear order on the 
set $X$ of generators of $L$, given by \fl(gen),
in the following way: 
\begin{equation}\label{fl:order}
a(m) < b(n) \ \iff   m<n \ \text{ or } \ 
( m=n  \ \text{ and } \   a<b ).
\end{equation}
On the set $X^*$ of words in $X$ introduce 
the standard lexicographical order: for $u,v \in X^*$ 
if $|u| < |v|$, set $u<v$; if 
$|u| = |v|$  then set $u<v$ whenever there is some $1\le i \le |v|$
such that $u(i) < v(i)$ and $u(j) = v(j)$ for all $1\le j < i$.

In a  defining relation from \fl(rel) the biggest term has form 
$b(n)a(m)$ such that 
\begin{equation}\label{fl:appl}
\text{$n-m > N$ or ($n-m = N$  and  
($b>a$  or  ($b=a$  and $N$ is odd)))}.
\end{equation}
Taking it as a principal part we get a rule on $\Bbbk\<X\>$:
\begin{subequations}\label{fl:therule}
\begin{equation}\label{fl:therulea}
\rho(b(n), a(m))  = 
\(b(n) a(m),\  
a(m) b(n) -  \sum_{s=1}^N  (-1)^s \binom Ns \ad{b(n-s)}{a(m+s)}\),
\end{equation}
and in case when $a=b, \, n-m=N$ and $N$ is odd, 
\begin{equation}\label{fl:theruleb}
\begin{split}
&\rho(a(m+N), a(m))  = \\
&\qquad\(a(m+N) a(m),\  
a(m) a(m+N) -  \frac 12 \sum_{s=1}^{(N-1)/2}  
(-1)^s \binom Ns \ad{a(n-s)}{a(m+s)}\).
\end{split}
\end{equation}
\end{subequations}

Denote the set of all such rules by $\mathcal R$:
\begin{equation}\label{fl:R}
\mathcal R = \{ \rho(b(n), a(m)) \ | \ \text{\fl(appl) holds}\}.
\end{equation}

\begin{Lem}\label{lem:confluency}\sl
The set of rules $\mathcal R$ is a confluent rewriting system on 
$\Bbbk\<X\>$.
\end{Lem}

We prove this Lemma in \sec{theproof}. Here we derive from it and
from the Diamond \lem{diamond} the following theorem.

\begin{Thm}\label{thm:basisass}\sl
\begin{itemize}
\item[(a)]
Let $C=C(N)$ be the free Lie conformal algebra generated by a linearly
ordered set
$\mathcal B$ corresponding to a constant locality function $N$. 
Let $L = \op{Coeff}C$ be the Lie algebra of coefficients, 
and let $U = U(L)$ be its 
universal enveloping algebra. Then a basis of $U$ is given by 
all monomials
\begin{equation}\label{fl:basisass}
a_1(n_1) a_2(n_2) \cdots a_k(n_k),\quad a_i \in \mathcal B, \ 
n_i \in \Z,
\end{equation}
such that for any $1\le i < k$ one has
\begin{equation}\label{fl:nojumps}
n_i - n_{i+1} \le 
\begin{cases}
N-1 & \text{if $a_i > a_{i+1}$ or ($a_i = a_{i+1}$ and $N$ is odd)}\\
N &\text{otherwise}.
\end{cases}
\end{equation}
\item[(b)]
A basis of the algebra $U(L_+)$ is given by all monomials \fl(basisass)
satisfying the condition \fl(nojumps) and such that all $n_i \ge 0$.
\item[(c)]
Let $V = V(C)$ be the corresponding free vertex algebra. Then a
basis of $V$ consists of elements 
\begin{equation}\label{fl:basisver}
a_1(n_1) a_2(n_2) \cdots a_k(n_k)\1,\quad a_i \in \mathcal B, \ 
n_i \in \Z,
\end{equation}
such that the condition \fl(nojumps) holds and, in addition, $n_k < 0$.
\end{itemize}
\end{Thm}

\begin{proof}
The statement
(a) is a direct corollary of \lem{confluency} and the Diamond \lem{diamond},
because \fl(basisass) is precisely the set of all terminal monomials
with respect to $\mathcal R$.

(b) follows immediately from \lem{Aplus}, since any subset of rules
$\mathcal R$ is also a confluent rewriting system. Note also that 
for a rule $\rho$ given by \fl(therule) if the principal term 
$\bar \rho$ contains only elements from $X_0$ then so does the whole
rule $\rho$.

For the proof of (c) recall  that $V \cong U/UL_+$ as linear spaces 
(and even as $L$-modules), where $UL_+$ is the
left ideal generated by $L_+$, see \sec{univertex}. By \lem{Aplus}, 
this ideal is the linear span of all
monomials $a_1(n_1) a_2(n_2) \cdots a_k(n_k)$ such that 
$n_k \ge 0$. But under the action of the rewriting system $\mathcal R$
the index of the rightmost symbol in a word can only increase, 
hence the linear span of these monomials in $\Bbbk\<X\>$ is
stable under $\mathcal R$.
It follows that the terminal monomials with a
non-negative rightmost index form a basis of $UL_+$. This
proves (b).
\end{proof}

\bigskip
\subsection{Proof of \lem{confluency}}
\label{sec:theproof}
First we prove that the set of rules  $\mathcal R$, given by \fl(R), is 
a rewriting system on $\Bbbk\<X\>$. Take a word 
$u = a_1(m_1)\cdots a_k(m_k) \in X^*$. Let $p\in \Bbbk\<X\>$  be such that
$u \ra p$. Then any word $v$ that appears in $p$ lies in the finite set
\begin{equation}\label{fl:Wu}
W_u = \left\{ b_1(n_1)\cdots b_k(n_k)\in X^*\ \left|\ 
n_i \ge \min_{1\le j\le k} \{m_j\} \text{  and  }
 \sum n_i = \sum m_i    
\right.\right\},
\end{equation}

therefore the condition (a) of the \dfn{rewriting} holds.

Thus we are left to prove that $\mathcal R$ is confluent. 
According to the Diamond \lem{diamond}, it is enough to check that
it is confluent on a composition $w=c(k)b(j)a(i)$ of principal parts
of a pair of rules $\rho(b(j),a(i)), \rho(c(k),b(j)) \in \mathcal R$. Thus
it is sufficient to prove the following claim.
\begin{Lem}\label{lem:confind}\sl
Let $u = c(k)b(j)a(i)\in X^*$ be a word of length 3. Then 
$\mathcal R$ is confluent on $u$, i.e. there is the unique
terminal $r(w)\in \Bbbk\<X\>$ such that $u \ra r(w)$.
\end{Lem}

\begin{proof}
Assume for simplicity that the three rules $\rho(b(n), a(m)), \ 
\rho(c(p), b(n))$ and $\rho(c(p), a(m))$ are of the form 
\fl(therulea). The general case is essentially the same, 
but requires some additional calculations.

Consider  the set $W_u$, given by \fl(Wu). 
We prove that the Lemma holds for all $w \in W_u$ by induction on $w$. If
$w$ is sufficiently small 
then it is a terminal itself. By induction, it is enough to 
consider $w = c(p)b(n)a(m) \in W_u$ such that $\mathcal R$ is applicable
to both $b(n)a(m)$ and $c(p)b(n)$. Apply $\rho(b(n), a(m))$ and 
$\rho(c(p), b(n))$ to $w$ and take the difference of the results: 
\begin{equation*}
\begin{split}
v = &\, b(n)c(p)a(m) - 
\sum_{s=1}^N (-1)^s \binom Ns \ad{c(p-s)}{b(n+s)} a(m) \\
 &- c(p)  a(m)b(n) + 
\sum_{s=1}^N (-1)^s \binom Ns c(p)\ad{b(n-s)}{a(m+s)}.
\end{split}
\end{equation*}

By induction,
$v$ is reduced uniquely to a terminal $t$ and we only have to show  
that $t=0$.

First we  apply the rules   
$\rho(b(n), a(m)),\ \rho(c(p), b(n))$ and $\rho(c(p), a(m))$ 
to $v$ several times and get
\begin{equation}\label{fl:to_prove}
\begin{split}
v \ra
- &\sum_{s=1}^N (-1)^s \binom Ns b(n)\ad{c(p-s)}{a(m+s)} 
+ b(n)a(m)c(p)\\
- &\sum_{s=1}^N (-1)^s \binom Ns \ad{c(p-s)}{b(n+s)} a(m) \\
+ &\sum_{s=1}^N (-1)^s \binom Ns \ad{c(p-s)}{a(m+s)} b(n)
- a(m)c(p)b(n)\\ 
+ &\sum_{s=1}^N (-1)^s \binom Ns c(p)\ad{b(n-s)}{a(m+s)}\\
\ra &\sum_{s=1}^N  (-1)^s \binom Ns \ad{a(m)}{\ad{c(p-s)}{b(n+s)}} \\
+ &\sum_{s=1}^N  (-1)^s \binom Ns \ad{\ad{c(p-s)}{a(m+s)}}{b(n)} \\
+ &\sum_{s=1}^N  (-1)^s \binom Ns \ad{c(p)}{\ad{b(n-s)}{a(m+s)}} 
\end{split}
\end{equation}

\def\jac{\kappa}
\def\loc{\lambda}

Next we introduce two rules acting on the linear combinations of
(formal) commutators: For any 
$ a(m), b(n), c(p) \in X $ let 
\begin{equation*}
\begin{split}
\jac &= \Big( \ad{a(m)}{\ad{b(n)}{c(p)}}, \ \ 
\ad{\ad{a(m)}{b(n)}}{c(p)} + \ad{b(n)}{\ad{a(m)}{c(p)}} \Big) \\
\loc &= \Big( \ad{b(n)}{a(m)},\ \ 
- \sum_{s=1}^N  (-1)^s \binom Ns \ad{b(n-s)}{a(m+s)} \Big)
\end{split}
\end{equation*}
The rule $\loc$ is the locality relation, and $\jac$ is nothing else 
but the Jacoby identity. The Lemma will be proved after we show two
things:
\begin{itemize}
\item[1)]\sl
There always exists a finite sequence of applications of the rules
$\jac$ and $\loc$ that reduces \fl(to_prove) to 0.
\item[2)]\sl
All words which appear in the process of reduction in {\rm 1)} are 
smaller than the initial word $u = c(p)b(n)a(m)$ with respect to the
order \fl(order).
\end{itemize}

Indeed, assume 1) and 2) hold. 
Denote the polynomial in \fl(to_prove) by $p_0$. Let 
$$
p_0 \ra p_1 \ra \dots \ra 0
$$
be the reduction, guaranteed by 1). By 2) and by the induction
hypothesis, any two neighboring polynomials 
$p_i \ra p_{i+1}$ from this sequence are uniquely $\mathcal R$-reduced
to a terminal, and this terminal must be the same, since either 
$p_i \xrightarrow{\mathcal R} p_{i+1}$ or 
$p_{i+1} \xrightarrow{\mathcal R} p_i$.

Denote the three last terms in \fl(to_prove) by \fbox{a}\,,\, \fbox{b} and
\fbox{c}\,. In Figure 1 we present a scheme of how 
$\jac$ and $\loc$ should be applied in order to reduce \fl(to_prove) to 0. 

\def\da{{\Large $\downarrow$}}
\def\fb#1{\fbox{#1}}

\begin{figure}[h]\label{fig:scheme}
\begin{center}
\small
\begin{tabular}{@{\extracolsep{-5pt}}rcccccccccccc}
      &      &\fb{a}&      &  +   &      &\fb{b}&      &  +   &      &\fb{c}&      &      \\
$\jac$:\ \  &      & \da  &      &      &      & \da  &      &      &      & \da  &      &      \\
      &\fb{d}&  +   &\fb{e}&      &\fb{f}&  +   &\fb{g}&      &\fb{h}&  +   &\fb{i}&      \\
$\loc$:\ \   & \da  &      & \da  &      & \da  &      & \da  &      & \da  &      & \da  &      \\
      &\fb{j}&      &\fb{k}&      &\fb{l}&      &\fb{m}&      &\fb{n}&      &\fb{o}&      \\
$\jac$:\ \   &      &      & \da  &      &      &      & \da  &      &      &      & \da  &      \\
      &      &\fb{p}&  +   &\fb{q}&      &\fb{r}&  +   &\fb{s}&      &\fb{t}&  +   &\fb{u}\\
$\loc$:\ \   &      & \da  &      &      &      &      &      & \da  &      &      &      & \da  \\
      &      &\fb{v}&      &      &      &      &      &\fb{w}&      &      &      &\fb{x}\\
$\jac$:\ \   &      &      &      &      &      &      &      & \da  &      &      &      &      \\
      &      &      &      &      &      &      &\fb{y}&   +  &\fb{z}&      &      & 
\end{tabular}
\end{center}
\caption{Application of rules $\jac$ and $\loc$} 
\end{figure}

Each box stands for a sum of commutators:

\begin{gather*}
\fbox{j} = - \, \fbox{r} = \sum_{s,t=1}^N (-1)^{s+t}\binom Ns \binom Nt 
\ad{\ad{c(p-s-t)}{a(m+t)}}{b(n+s)},  \\
\fbox{k} = \sum_{s,t=1}^N (-1)^{s+t}\binom Ns \binom Nt 
\ad{c(p-s)}{\ad{b(n+s-t)}{a(m+t)}},  \\
\fbox{l} = - \, \fbox{t} = \sum_{s,t=1}^N (-1)^{s+t}\binom Ns \binom Nt 
\ad{c(p-s)}{\ad{b(n-t)}{a(m+s+t)}},  \\
\fbox{m} = \sum_{s,t=1}^N (-1)^{s+t}\binom Ns \binom Nt 
\ad{\ad{b(n+t)}{c(p-s-t)}}{a(m+s)},  \\
\fbox{n} = - \, \fbox{q} = \sum_{s,t=1}^N (-1)^{s+t}\binom Ns \binom Nt 
\ad{\ad{b(n-s+t)}{c(p-t)}}{a(m+s)},  \\
\fbox{o} = \sum_{s,t=1}^N (-1)^{s+t}\binom Ns \binom Nt 
\ad{b(n-s)}{\ad{a(m+s+t)}{c(p-t)}}, \\
\fbox{v} = - \, \fbox{y} = \sum_{s,t,r=1}^N (-1)^{s+t+r}\binom Ns \binom Nt 
\binom Nr \ad{b(n+s-t)}{\ad{a(m+t+r)}{c(p-s-r)}}, \\
\fbox{w} = \sum_{s,t,r=1}^N (-1)^{s+t+r}\binom Ns \binom Nt 
\binom Nr \ad{\ad{a(m+s+r)}{b(n+t-r)}}{c(p-s-t)}, \\
\fbox{x} = - \, \fbox{z} = \sum_{s,t,r=1}^N (-1)^{s+t+r}\binom Ns \binom Nt 
\binom Nr \ad{a(m+s+t)}{\ad{c(p-t-r)}{b(n-s+r)}}, \\
\end{gather*}

One can see that all terminal boxes in the above scheme cancel, so that
$\fbox{a} + \fbox{b}+ \fbox{c} \ra 0$. 
Claim 2) also holds, since every symbol in every box
in Figure 1 is less than $c(p)$. 

\end{proof}

\bigskip
\subsection{Digression on Hall bases}
\label{sec:Hall}
Let again $\mathcal B$ be some linearly ordered alphabet, $N\in \Z_+$, 
$C = C(N)$ the free Lie conformal algebra generated by $\mathcal B$
with respect to the constant locality $N$, and $L=\op{Coeff} C(N)$.
A basis of the Lie algebra $L$ could be obtained
by modifying the construction of a Hall basis of a free Lie algebra,
see \cite{hall}, \cite{shirshov58}, \cite{shirshov62a}.
Here we review the latter construction. We closely follow \cite{reut},
except that all the order relations are reversed.

As in \sec{diamond}, take an alphabet $X$ and a commutative
ring $K$. Let $T(X)$ be the set of all 
binary trees with 
leaves from  $X$. For typographical reasons we will write 
the tree $\^{x\ y}$ as $\tr xy$.  Assume that $T(X)$ is endowed
with a linear order such that  $\tr xy > \min\{x,y\}$ for any
$x,y \in T(X)$.

\begin{Def}\label{dfn:Hall}
A {\it Hall set} $\mathcal H \subset T(X)$ is a subset of all
trees $h \in T(X)$ satisfying the following (recursive) 
properties:
\begin{itemize}
\item[1.]
If $h = \tr xy$ then $y,x \in \mathcal H$ and $x>y$;
\item[2.]
If $h = \tr{\tr xy}z$ then $z \ge y$, \ so that 
$\tr xy > z \ge y$.
\end{itemize}
In particular, $X \subset \mathcal H$.
\end{Def}

Introduce two maps 
$\alpha: T(X) \to X^*$ and 
$\lambda : T(X) \to K\<X\>$
in the following recursive way: 
for $a \in X$ set $\alpha(a)=\lambda(a) = a$
and $\alpha(\tr xy) = \alpha(x)\alpha(y), \ 
\lambda(\tr xy) = \ad{\lambda(x)}{\lambda(y)}$. 

It is a well-known fact (see e.g. \cite{reut}) that 
\begin{itemize}\sl
\item[(a)]
$\lambda(\mathcal H)$ is a basis of the  free Lie algebra generated by 
$X$ and
\item[(b)]
$\alpha\big|\raisebox{-5pt}{$\mathcal H$}$ is injective.
\end{itemize}

A word $w \in \alpha(\mathcal H)$ is called a {\it Hall word}.

On the set $X^*$ of words in $X$ introduce a
(lexicographic) order as follows: if $u$ is a prefix of $v$ then 
$u > v$, otherwise $u>v$ whenever for some index $i$ one has 
$u_i > v_i$ and $u_j = v_j$ for all $j < i$. 

\begin{Def}\label{dfn:Lyndon} (\cite{shirshov62b}, \cite{lyndon58})
A word $v \in X^*$
is called {\it Lyndon-Shirshov} if it is bigger than all its proper suffices.
\end{Def}

\begin{Prop}\label{prop:Lyndon}\sl
\begin{itemize}
\item[(a)]
There is a Hall set $\mathcal H_{\rm LS}$ such that 
$\alpha(\mathcal H_{\rm LS})$ is the set of all Lyndon-Shirshov words
and $\alpha: T(X) \to X^*$ preserves the order.
\item[(b)]
For any tree $h \in \mathcal H_{\rm LS}$ the biggest term in 
$\lambda(h)$ is $\alpha(h)$.
\end{itemize}
\end{Prop}

\bigskip
\subsection{Basis of the algebra of coefficients of a free
Lie conformal algebra}
\label{sec:basislie} 
Here we apply general results from \sec{Hall} to the situation of 
\sec{freeconf}.

Recall that starting from a set of symbols $\mathcal B$ and a number 
$N> 0$, we build the free conformal algebra $C = C(N)$ generated
by $\mathcal B$ such that $a \encirc n b = 0$ for any two 
$a,b \in \mathcal B$ and $n\ge N$. Let  $L = \op{Coeff}C$ be the 
corresponding Lie algebra of coefficients. It is generated by the set
$X = \{a(n)\ | \ a\in \mathcal B, \ n \in \Z \}$ subject
to relations \fl(rel). 

The set of generators $X$ is equipped with the linear 
order defined  by \fl(order). We define the order on $X^*$
as in \sec{Hall}. Consider the set of all Lyndon words in 
$X^*$ and let 
$\mathcal H = \mathcal H_{\text{LS}} \subset T(X)$ be the 
corresponding Hall set. Recall that there is  a rewriting system 
$\mathcal R$ on $\Bbbk\<X\>$, given by \fl(R). Define
$$
\mathcal H_{\rm term} = 
\{h \in \mathcal H \ | \ \alpha(h) \text{  is terminal}\}.
$$

\begin{Lem}\label{lem:basislie}\sl
\begin{itemize}
\item[(a)]
Let $v_1 \le \dots\le v_n$ be a non-decreasing sequence of terminal
Lyndon-Shirshov words. Then their concatenation 
$w = v_1\cdots v_n \in X^*$ is a terminal word.
\item[(b)]
Each terminal word $w \in X^*$ can be uniquely represented 
as a concatenation $w = v_1\cdots v_n$, where $v_1 \le \dots\le v_n$
is a non-decreasing sequence of terminal Lyndon-Shirshov words.
\end{itemize}
\end{Lem}

\begin{proof}
(a)  Take two terminal Lyndon-Shirshov words $v_1 \le v_2$. Let $x \in X$
be the last symbol of $v_1$ and $y \in X$
be the first symbol of $v_2$. Then, since a word is less than its prefix
and since $v_1$ is a Lyndon-Shirshov word, we get 
$$
x < v_1 \le v_2 < y.
$$
Therefore, $xy$ is a terminal, hence $v_1v_2$ is a terminal too.

(b)
Take a terminal word $w \in X^*$. Assume it is not Lyndon-Shirshov.
Let $v$ be the maximal among
all proper  suffices of $w$. Then $v$ is Lyndon-Shirshov, $v > w$ and $w = uv$
for some word $u$. By induction,
$u = v_1 \dots v_{n-1}$ for a non-decreasing sequence of Lyndon-Shirshov words
$v_1 \le \dots\le v_{n-1}$. We are left to show that  $v \ge v_{n-1}$.

Assume on the contrary that  $v < v_{n-1}$. Then, since $v > v_{n-1}v$, 
\ $v_{n-1}$ must be a prefix of $v$ so that $v = v_{n-1} v'$. But then
$v' > v$ which contradicts the Lyndon-Shirshov property of $v$. 

The uniqueness is obvious.
\end{proof}

Let $\varphi : \Bbbk\<X\> \to U(L)$ be the canonical projection
with the kernel $I(\mathcal R)$.
\begin{Thm}\label{thm:basislie}\sl
The set $\varphi(\lambda(\mathcal H_{\rm term}))$ is a basis of $L$.
\end{Thm}

\begin{proof}
Let $s = \{h_1, \dots, h_n\} \subset \mathcal H_{\rm term}$ be 
a non-decreasing sequence of terminal Hall trees. Let 
$\lambda(s) = \lambda(h_1)\cdots\lambda(h_n) \in \Bbbk\<X\>$
and $\alpha(s) = \alpha(h_1)\cdots\alpha(h_n) \in X^*$

By the Poincar\'e-Birkoff-Witt theorem 
it is sufficient to prove that  the set
$\{\varphi(\lambda(s))\}$, when $s$ ranges over all non-decreasing
sequences $s$ of  terminal Hall trees, is a basis of $U(L)$. 

By (b) of \prop{Lyndon}, $\lambda(s) = \alpha(s) + O(\alpha(s))$,
where $O(v)$ stands for a sum of terms which are less than $v$.
Now let $t(s)\in \Bbbk\<X\>$ be a terminal such that 
$\lambda(s) \ra t(s)$. One can view $t(s)$ as the decomposition of 
$\varphi(\lambda(s))$ in basis \fl(basisass). By \lem{basislie},
$\alpha(s)$ is a terminal monomial, hence $t(s)$ has form
$t(s) = \alpha(s) + f(s)$ where $f(s)$ is a sum of terms 
$v\in X^*$
satisfying the following properties: 
\begin{itemize}
\item[1.]
$v$ is terminal and $v<\alpha(s)$;
\item[2.]
If $v$ contains a symbol $a(n) \in X$ then
$a$ appears in $\alpha(s)$ and $n_{\min} \le n \le n_{\max}$, 
where $n_{\min}$ and $n_{\max}$ are respectively minimum and maximum
of all indices that appear in $\alpha(s)$.
\end{itemize}
Indeed, due to \prop{Lyndon} b) properties 1 and 2 are satisfied by 
all the terms in $\lambda(s) - \alpha(s)$ , and they cannot be broken 
by an application of the rules $\mathcal R$.

Property 1 implies that all $t(s)$ and,
therefore, $\varphi(\lambda(s))$ are linearly independent.
So we are left to show that they span $U(L)$. For that purpose 
we show that any terminal word $w \in X^*$ can be represented
as a linear combination of $t(s)$. 

By (b) of \lem{basislie} any terminal word $w$ could be written 
as $w = \alpha(s)$ for some non-decreasing sequence $s$ of terminal
Hall trees. So we can write $w = t(s) - f(s)$. Now do the same 
with any term $v$ that appear in $f(s)$, and so on. This process
should terminate, because every term $v$ that appears during this process
must satisfy properties 1 and 2 and there are only finitely many such
terms.

\end{proof}

\begin{Rem}
Altrernatively we could use the theorem of L. Bokut' and P. Malcolmson
\cite{BM}.
\end{Rem}

As in (b) of \thm{basisass}, we deduce that all the elements of 
$\varphi(\lambda(\mathcal H_{\rm term}))$ containing only symbols
from $X_0$ form a basis of $L_+$.

Note that we have an algorithm for building a basis of the 
free Lie conformal algebra $C = C(N)$. Let $L = \op{Coeff}C$, 
$V = V(C)$ and $U = U(L)$. Recall that the image if $C$ in 
$V$ under the canonical embedding $\psi: C\to V$ is 
$\psi(C) = L_-\1 = L\1 \subset V$.
So, the algorithm goes as follows: take the basis of $L$ 
provided by \thm{basislie}. Decompose its element in basis
\fl(basisass) of the universal enveloping algebra $U(L)$, and then
cancel all terms of the form $a_1(n_1)\cdots a_k(n_k)$ where $n_k \ge 0$.  
What remains, being interpreted as elements of the vertex algebra
$V$, form a basis of $\psi(C)\subset V$.

\bigskip
\subsection{Basis of the algebra of coefficients of a free
associative conformal algebra}
\label{sec:basisass}
Let again $\mathcal B$ be some alphabet, and 
$N:\mathcal B \times \mathcal B \to \Z_+$ be a locality
function, not necessarily constant and not necessarily symmetric.
By \prop{freeconf}, the coefficient algebra $A = \op{Coeff} C(N)$
of the free associative conformal algebra $C(N)$ corresponding
to the locality function $N$ is presented in terms of generators and
relations by the set of generators 
$X = \{ b(n) \ | \ b\in \mathcal B, n \in \Z \}$
and relations \fl(rel).

\begin{Thm}\label{thm:assbasis}\sl
\begin{itemize}
\item[(a)]
A basis of the algebra $A$ is given by all monomials of the form
\begin{equation}\label{fl:basisa}
a_1(n_1)\cdots a_{l-1}(n_{l-1})a_l(n_l), 
\end{equation}
where $a_i \in \mathcal B$ and 
$$
-\left\lceil\frac{N_i-1}2\right\rceil \le n_i 
\le \left\lfloor\frac{N_i-1}2\right\rfloor, \ 
N_i = N(a_i, a_{i+1}), \ \text{  for  }
i= 1,\dots, l-1.
$$
\item[(b)]
A basis of the algebra $A_+$ is given by all monomials of the form
\begin{equation}\label{fl:basisaplus}
a_1(n_1)\cdots a_{l-1}(n_{l-1})a_l(n_l), 
\end{equation}
where $a_i \in \mathcal B$ and 
$$
0\le n_i \le N_i-1, \ 
N_i = N(a_i, a_{i+1}), \ \text{  for  }
i= 1,\dots, l-1.
$$
\end{itemize}
\end{Thm}

\begin{Cor}\label{cor:assbasis}\sl 
Assume that the locality function $N$ is constant. Consider
the homogeneous component $A_{k,l}$ of $A$, spanned
by all the words of the length $l$ and of the sum of indexes 
$k$. Then $\dim A_{k,l} = N^{l-1}$.
\end{Cor}

\begin{proof}[Proof of \thm{assbasis}]
(a)\quad
Introduce a linear order on $\mathcal B$, and define an order on the set
of generators $X$ by the following rule:
$$
a(m) > b(n) \iff |m| > |n| \text{ or }
m=-n>0 \text{ or }
(m=n  \text{ and } a>b)
$$
In particular, for some $a\in \mathcal B$ we have
$$
a(0) < a(-1) < a(1) < a(-2) < a(2) < \ldots.
$$

For any relation $r$ from \fl(rel) take the biggest term 
$\bar r$ and consider the rule $(\bar r,\ r-\bar r)$. This way
we get a collection of rules 
\begin{equation*}
\begin{split}
\mathcal R  = &\left\{ \rho_1(b(n), a(m)) \ \bigg| \ 
a,b \in \mathcal B, \ 
n >  \left\lfloor\tfrac{N(b,a)-1}2\right\rfloor \right\} \bigcup \\
&\left\{ \rho_2(b(n), a(m)) \ \bigg| \ 
n < -\left\lceil\tfrac{N(b,a)-1}2\right\rceil \right\}
\end{split}
\end{equation*}
where
\begin{gather*}
\rho_1(b(n), a(m)) = 
\left( b(n)a(m), \ \  \sum_{s=1}^{N(b,a)} (-1)^{s+1}  
\binom{N(b,a)}{s} b(n-s)a(m+s) \right), \\
\rho_2(b(n), a(m)) = 
\left( b(n)a(m), \ \  \sum_{s=1}^{N(b,a)} (-1)^{s+1}  
\binom{N(b,a)}{s} b(n+s)a(m-s) \right).
\end{gather*}

By the Diamond \lem{diamond}, 
we have to prove that these rules form a confluent rewriting
system on $\Bbbk\<X\>$. 
Clearly $\mathcal R$ is a rewriting system, since it decreases 
the order, and each subset of $\Bbbk\<X\>$, containing only finitely
many different letters from $\mathcal B$, has the minimal
element, in contrast to the situation of \sec{theproof}.

As before, it is enough to check that $\mathcal R$
is confluent on 
any composition $w = c(p)b(n)a(m)$,
of the principal parts of rules from $\mathcal R$. 
Consider the set 
$W = \{ \ c(k)b(j)a(i) \ | \ k,j,i \in \Z\ \}\subset X^*$.
we prove by induction on $w\in W$ that 
$\mathcal R$ is confluent on $w$. If $w$ is sufficiently small, then
it is terminal. Assume that $w = c(k)b(j)a(i)$ is an ambiguity, 
for example that $\rho_1(c(p), b(n))$ and $\rho_2(b(n), a(m))$
are both applicable to $w$. Other cases are done in the same way. 
Let 
\begin{gather*}
w_1 = \rho_1(c(p), b(n))(w) = 
\sum_{s=1}^{N(c,b)} (-1)^s \binom{N(c,b)}{s} c(p-s)b(n+s)a(m), \\
w_2 = \rho_2(b(n), a(m))(w) = 
\sum_{t=1}^{N(b,a)} (-1)^t \binom{N(b,a)}{t} c(p)b(n+t)a(m-t).
\end{gather*}
Applying $\rho_2(b(n+s), a(m))$  for $s=1, \ldots, N(b,a)$ to $w_1$
gives the same result as we get from 
applying $\rho_1(c(p), b(n+t))$  for $t=1, \ldots, N(c,b)$ to $w_2$, 
namely
\begin{equation}\label{fl:common}
\sum_{s,t\ge1}(-1)^{s+t} \binom{N(c,b)}{s}\binom{N(b,a)}{t} 
c(p-s)b(n+s+t)a(m-t).
\end{equation}
By the induction assumption, $w_1 - w_2$ is uniquely reduced to a
terminal, and since all monomials in \fl(common) are smaller than 
$w$, we conclude that this terminal must be 0.

\smallskip
\noindent(b)\quad
Follows at once from \lem{Aplus}.
\end{proof}

\bigskip
\section*{Acknowledgments}
I am grateful to Bong Lian, Efim Zelmanov and Gregg Zuckerman
for helpful discussions. I thank Victor Kac, Bong Lian and Gregg 
Zuckerman for communicating unpublished papers (\cite{kac_loc}, 
\cite{kac_fd}, \cite{lz_preprint}) and I thank L.~A.~Bokut' who
carefully read this paper and made valuable comments.

\bigskip
\bibliography{../conformal}

\end{document}